\documentclass[11pt]{article}
\usepackage{amssymb,latexsym,amsmath,amsbsy,amsthm,amsxtra,amsgen,dsfont,pdfsync,color}
\begin{document}
\title{On the genealogical structure of \\ critical branching processes in a varying environment}

\author{G\"otz Kersting\thanks{Institut f\"ur Mathematik, Goethe Universit\"at, Frankfurt am Main, Germany, kersting@math.uni-frankfurt.de}}

\maketitle
\begin{abstract}
Critical branching processes in a varying environment behave much the same as critical Galton-Watson processes. In this note we like to confirm this finding with regard to the underlying genealogical structures. In particular, we consider the most recent common ancestor given survival and the corresponding reduced branching processes, in the spirit of Zubkov (1975) and Fleischmann and Siegmund-Schultze (1977). 
\medskip

\noindent
\textit{Keywords and phrases.}  branching process, varying environment, Galton-Watson process, critical process, reduced process, most recent common ancestor,  exponential distribution, Yule process

\smallskip
\noindent
\textit{MSC 2010 subject classification.} Primary  60J80.\\
\end{abstract}

\section{Introduction}
Critical Galton-Watson processes (GW processes)  constitute  a prominent, but small portion of the whole family of GW processes, and could be considered just as the separating borderline between the larger classes of supercritical and subcritical GW processes. From this point of view it might come as a surprise that within the vaster family of branching processes in a varying environment the subclass of critical processes is much broader and  pushes apart the super- and subcritical processes. Here we consider a branching process in a varying environment to be critical if it shares properties characteristic of critical GW processes. As in \cite{Ke}, we think of Kolmogorov's asymptotics on the survival probabilities and Yaglom's asymptotic exponential law given survival. In this note we like to back this approach by widening it  to the genealogical structure  of the underlying branching trees, in the spirit of Zubkov's result \cite{Zu} on the most recent common ancestor within a critical  GW tree and the  treatment of its reduced process by Fleischmann and Siegmund-Schultze \cite{FlSi}. Here again critical branching processes in a varying environment resemble critical GW processes.

In order to define a branching process in a varying environment (BPVE), let $Y_1, Y_2, \ldots$ denote a sequence of random variables with values in $\mathbb N_0$, and  $f_1,f_2, \ldots$ their distributions. Let $Y_{in}$, $i,n\in \mathbb N$, be independent random variables such that $Y_{in}$ and $Y_n$ coincide in distribution for all $i,n\ge 1$. Define the random variables $Z_n$, $n \ge 0$, with values in $\mathbb N_0$ recursively as
\[ Z_0 :=1 \ , \quad Z_{n} := \sum_{i=1}^{Z_{n-1}} Y_{in} \ , \ n \ge 1\ . \]
Then the process $(Z_n)_{n \ge 0}$ is called a {\em branching process in the varying environment $v=(f_1,f_2, \ldots)$ with initial value $Z_0=1$.} These processes may be considered as a model for the development of the size of a population where individuals reproduce independently  with offspring distributions $f_n$ potentially changing among generations. Without further mention we always require that  $0 < \mathbf E[Y_n] < \infty$ for all $n \ge 1$. 

In this introduction we recall the notion of a critical BPVE, as it was initiated by Jagers in \cite{Ja} and  developed only recently  both by Bhattacharya and Perlman \cite{BhaPe} and by Kersting \cite{Ke}. First, we separate certain exceptional BPVEs out,  behaving in an extraordinary manner  unknown for ordinary GW processes (like possessing multiple rates of growth, see \cite{Schuh}). As explained in \cite{Ke}, these exotic possibilities can be precluded by certain uniform $L^2$-integrability conditions. Here we impose the requirement that for every $\varepsilon >0$ there is a constant $c_\varepsilon < \infty$ such that for all $n \ge 1$
\begin{align} 
\mathbf E\big [ Y_n^2 ; Y_n > c_\varepsilon(1+ \mathbf E[Y_n])\big] \le \varepsilon \mathbf E \big[ Y_n^2; Y_n \ge 2\big] < \infty \ ,
\tag{$\ast$} 
\end{align}
which is introduced in \cite{Ke} as condition (B). There it is demonstrated in detail that this assumption is widely satisfied. Let us mention  just two examples: if either the $Y_n$ are uniformly bounded by some constant $c<\infty$, or if else all $Y_n$ are Poisson distributed (without any constrain on the parameters), then ($\ast$) holds true. In particular,  the latter example   shows that condition ($\ast$) does not entail any restriction on the sequence of expectations of $Y_n$ respectively $Z_n$.

It turns out that under assumption ($\ast$) the behaviour of a BPVE $(Z_n)_{n \ge 0}$ is essentially dictated by two sequences of numbers, namely
\[ \mu_n:= \mathbf E[Z_n] \ \text{ and }   \ \rho_n := \frac{\mathbf E[Z_n(Z_n-1)]}{\mathbf E[Z_n]^2} \ , \ n \ge 0\  .  \]
With these quantities the family of BPVEs can be classified into subcritical, critical, supercritical and asymptotically degenerate processes (where the process may freeze in a positive state). For complete details we refer to \cite{Ke},  let us  recall just that statement  which is relevant in our context. It is a version of  \mbox{\cite[Theorem 4]{Ke}} (also taking \cite[Lemma 4]{Ke} into account).

\paragraph{Theorem 1.} {\em Assume} ($\ast$). {\em If the probability of extinction is 1, that is if $\mathbf P(Z_n>0) \to 0$ as $n \to \infty$, then the following statements are equivalent:
\begin{enumerate}
\item[\em (i)] There is a sequence $b_n$, $n \ge 0$, of positive numbers such that   $Z_n/b_n$ conditioned on the event $\{Z_n>0\}$ converges in distribution to a standard exponential distribution as $n \to \infty$,
\item[\em (ii)] $\mathbf E[Z_n \mid Z_n >0] \to \infty$ as $n \to \infty$,
\item[\em (iii)] $\mu_n\rho_n \to \infty$ as $n \to \infty$.
\end{enumerate}

\noindent
Under these conditions we may set $b_n:= \mathbf E[Z_n \mid Z_n >0]$, and we have
\[ \mathbf E[Z_n \mid Z_n >0] \sim   \frac 12 \mu_n \rho_n \ ,\]
or equivalently
\[ \mathbf P(Z_n>0) \sim \frac 2{\rho_n}  \]
as $n \to \infty$.}

\mbox{}\\
This theorem characterizes in the context of BPVEs the range, where Yaglom's \cite{Ya} asymptotics for critical GW processes on the conditional distribution of $Z_n$ holds. For a probabilistic approach via a 2-spine decomposition see \cite{PaTo}. In addition, the  final statement of the theorem provides an extension to BPVEs of Kolmogorov's \cite{Ko} asymptotics on the survival probability of  critical GW processes (in its general version, as given in \cite{KeNeSp}). Therefore, assuming ($\ast$) it is natural to transfere by means of Theorem 1 the notion of criticality  from GW processes to BPVEs, that is, to name a BPVE $(Z_n)_{n \ge 0}$ {\em critical} if both
\[ \mathbf P(Z_n>0) \to 0 \ \text{ and }\ \mathbf E[Z_n \mid Z_n >0] \to \infty  \]
as $n\to \infty$ are satisfied.  \cite[Proposition 1]{Ke} states that this pair of requirements can equivalently be expressed by the conditions
\begin{align} \rho_n \to \infty \ \text{ and }\ \mu_n\rho_n \to \infty  
\label{critical}
\end{align}
as $n \to \infty$. If the first one fails then we enter the region of supercritical or asymptotic degenerate processes, and if the second one goes wrong then we access the subcritical domain.

To conclude this introduction we recall how $\mu_n$ and $\rho_n$ may be derived from the varying environment $(f_1,f_2, \ldots)$.
For this purpose let us agree on the  notational facilitation used in \cite{Ke}: for a probability measure $f$ on $\mathbb N_0$ with weights   $f[k]$, $k \in \mathbb N_0$, we set
\[ f(s):= \sum_{k=0}^\infty s^k f[k] \ , \ 0 \le s \le 1\ . \]
Thus, we denote the probability measure $f$ and its generating function by one and the same symbol. Then, each operation applied to these measures has to be understood as an operation applied to their generating functions. For example,  $f_1 \circ f_2$ expresses the composition of  generating functions as well as the resulting probability measure. We shall consider the  mean and  the normalized second factorial moment  of a random variable $Y$ with distribution $f$,
\[ \mathbf E[Y]= f'(1)  \  \text{ and } \ \nu:= \frac{\mathbf E[Y(Y-1)]}{\mathbf E[Y]^2}= \frac{f''(1)}{f'(1)^2}\ . \]

With these convention in mind and setting
\[ \nu_n:= \frac{f_n''(1)}{f_n'(1)^2} \ , \]
we  obtain $\mu_n$ and $\rho_n$ by means of the formulas
\begin{align} \mu_n = f_1'(1)\cdots f_n'(1) \ , \quad \rho_n= \sum_{k=1}^n \frac{\nu_k}{\mu_{k-1}}\ . 
\label{murho}
\end{align}
Thus $(\rho_n)_{n \ge 0}$ is an increasing sequence, and also $\mu_0=1$ and $\rho_0=0$. For the proof we refer to \mbox{\cite[Lemma 4]{Ke}}.

\paragraph{Examples.} 1. Let $f_n$, $n \ge 1$, be a sequence of Poisson distributions with parameters $\lambda_n$. Then $f_n'(1)=\lambda_n$ and $\nu_n=1/\lambda_n$. Therefore, because of \eqref{critical} the BPVE is critical, if and only if
\begin{align}\sum_{k=1}^\infty \frac 1{\mu_k} =\infty  \ \text{ and } \ \frac 1{\mu_n} = o \Big(\sum_{k=1}^n \frac 1{\mu_k} \Big) \ . 
\label{Poissoncriterion}
\end{align}
The first condition is violated, if $\mu_n$ increases slightly faster than linearly, then we observe  supercritical behaviour. The second conditions fails, if $\mu_n$ decreases at an exponential rate, then we encounter subcritical processes. In between, a variety of behaviour for $\mu_n=\mathbf E[Z_n]$  can be set up. This demonstrates the broad range of critical BPVEs.

2. In the case of binary offspring  (meaning that $f_n[0]+ f_n[2]=1$ for all $n\ge 1$) we have once again $\nu_n=1/f_n'(1)$. Therefore, the  criterion \eqref{Poissoncriterion} for criticality applies equally. 

3. These examples are in a way typical and generalize as follows: Assume that for some $a>0$ we have for all $ n \ge 1$
\[ f_n''(1) \ge a f_n'(1) \ . \]
Additionally, from formula \eqref{secondfirst} below we see that ($\ast$) implies $f_n''(1) \le b f_n'(1)(1+f_n'(1))$ for some $b>0$. These estimates imply
\[ a\sum_{k=1}^n \frac 1{\mu_k} \le \sum_{k=1}^n \frac{\nu_k}{\mu_{k-1}}  \le b \sum_{k=1}^n \Big( \frac 1{\mu_k} + \frac 1{\mu_{k-1} }\Big)\le 2b \sum_{k=0}^n \frac 1{\mu_k} \ .\]
Then \eqref{critical} and \eqref{Poissoncriterion} are equivalent conditions.

\section{The genealogical structure of a critical BPVE}

Conditional  on the event that $Z_n >0$ we may define the generation $G_n$ of the most recent common ancestor of the $Z_n$ particles in generation $n$. Zubkov \cite{Zu} proved that for critical GW processes $G_n/n$ has asymptotically a uniform distribution on the interval $[0,1]$ (for another  proof see \cite{Ge}). The following theorem generalizes this remarkable result to  BPVEs.

\paragraph{Theorem 2.} {\em Assuming} ($\ast$) {\em we have for a critical BPVE
\[ \max_{0\le k \le n} \Big|\mathbf P(G_n \le k\mid Z_n >0) - \frac{\rho_k}{\rho_n}\Big| \to 0 \]
as $n \to \infty$.}

\mbox{}\\
Remarkably, the approximating cumulative distribution function $\rho_k/\rho_n$, $0\le k \le n$, also appears in
the recent paper \cite{PaTo} by Natalia Cardona-Tob\'on and Sandra Palau. There it determines  the moment of bifurcation of spines in a 2-spine construction of a (twofold) size-biased critical BPVE. One may wonder, whether Theorem 2 can be derived as well within this  framework.  

For a corollary we note that in \eqref{murho} the summands of $\rho_n$ are uniformly of smaller order than $\rho_n$ itself. This is the content of Lemma 1 below. Therefore $\mathbf P(G_n \le k \mid Z_n >0)\to 0$ as well as $\mathbf P(G_n \ge n-k \mid Z_n>0) \to 0$ for any $k \ge 1$ as $n \to \infty$. In other terms the following holds.

\paragraph{Corollary.} {\em Under the assumptions of Theorem 2, given the events $\{Z_n>0\}$ we have
\[ G_n \to \infty \ \text{ and } \ n-G_n \to \infty \]
in probability as $n \to \infty$.}

\mbox{}\\
These properties are characteristic:  For supercritical (and asymptotically degenerate) BPVRs one  expects that the random variables $G_n$, given that $Z_n>0$, remain  bounded in probability uniformly in $n$, and for subcritical the same will hold true for the quantities $n-G_n$.

The next examples discuss cases where, after a suitable normalization, $G_n$ given $\{Z_n>0\}$ has   a limiting distribution.

\paragraph{Examples.} 1. $\rho_{[tn]}/\rho_n$ converges for all $0<t<1$ to a limiting value, if and only if the sequence $(\rho_n)_{n\ge 1}$  is regularly varying with an exponent $\lambda \ge 0$. This means that for any $0<t<1$
\[ \mathbf P(G_n \le nt \mid Z_n >0) \to t^\lambda \]
as $n \to \infty$. Thus $G_n/n$ has a non-degenerate conditional limiting distribution, if and only if $(\rho_n)_{n\ge 1}$  is regularly varying with an exponent $\lambda >0$.

2. In the Poisson case, suppose that $\mu_n \sim n$, and so $\rho_n \sim \sum_{k=1}^n k^{-1} \sim \log n$. Thus \eqref{critical} is satisfied, and we are in a critical region close to supercritical behaviour. Then for $0<t<1$
\[ \mathbf P(\log G_n \le t \log n\mid Z_n>0) = \mathbf P(G_n \le n^t \mid Z_n>0) \sim \frac{\log n^t}{\log n}=t \ . \]
Here, the conditional distribution of $(\log G_n)/\log n$  is asymptotically uniform. 

3. In the Poisson case, suppose that $\mu_n \sim \exp(- \sqrt n)$. Since by l'Hospital's rule 
\[ \int_0^n \exp(\sqrt y)\, dy \sim 2 \sqrt n \exp(\sqrt n)  \]
as $n \to \infty$, we have $\rho_n \sim  2 \sqrt n \exp(\sqrt n)$. Again, we are in the critical domain, now close to  subcritical processes. In this case for $t>0$
\[ \mathbf P( G_n \le n- 2t\sqrt n \mid Z_n>0) \sim \frac {\rho_{n-[2t\sqrt n]}}{\rho_n} \sim \exp \Big(\sqrt{n-2t\sqrt n}-\sqrt n\Big) \sim \exp(-t) \ .\]
In other words: the conditional distribution of $(n-G_n)/(2\sqrt n)$ is asymptotically standard exponential. 

\mbox{}\\
From the  generation of the most recent common ancestor, we step forward to the more ela\-borate {\em reduced branching trees} given the events $\{Z_n >0\}$. As introduced by Fleischmann and Siegmund-Schultze \cite{FlSi} in the case of GW trees, these are the subtrees  containing all individuals which possess  descendants in generation $n$. Accordingly, the {\em reduced branching process} $(Z_{k,n})_{0\le k \le n}$ is made up of  the numbers $Z_{k,n}$ of individuals in generation $k\le n$ having  descendants in generation~$n$. Just as in the GW case, this process is a non-homogeneous Markov chain, actually a BPVE. For critical GW processes the  structure of the reduced process has been investigated in \cite{FlSi}. It turns out that asymptotically  the process $(Z_{[nt],n})_{0\le t < 1}$   coincides  with a properly time-changed Yule process. It is noteworthy that in  the case of a critical branching processes in a {\em random} environment  the situation is  markedly different, as clarified by Borovkov, \mbox{Vatutin \cite{BoVa}} and Vatutin \cite{Va}.

By contrast, for critical BPVEs we encounter  much the same behavior as for GW processes.  Again, we regard the increasing sequence $\rho_k/\rho_n$, $0 \le k \le n$, as a cumulative probability distribution function. We use its (right-continuous) generalized inverse \mbox{$k_n:[0,1)\to \{0,\ldots,n-1\}$} given by
\[ k_n(t):= \max \{ k \ge 0: \rho_k \le t \rho_n\} \ , \quad 0\le t<1 \ . \] 
By means of this function we rescale distances between generations in the reduced process. The following theorem contains the result of Fleischmann and Siegmund-Schultze as a special case.

\paragraph{Theorem 3.} {\em Assume} ($\ast$). {\em Then we have for the reduced process of a critical BPVE in the limit $n \to \infty$
\[    \Big( (Z_{k_n(t),n})_{0 \le t < 1} \, \big| \, Z_n>0 \Big) \ \stackrel{d} \Rightarrow \ \big(Y(\log \tfrac 1{1-t})\big)_{0\le t < 1 }\ ,   \]
where $(Y(u))_{u \ge 0}$ denotes a standard Yule process.}

\mbox{}\\
Here the notation $\stackrel d\Rightarrow$ indicates that for any $0<\eta < 1$ we have the convergence in the Skorohod sense of the sequence of processes restricted to times   $ t \in [0, \eta]$.

\section{Proofs}
We denote by $f_{k,n}$, $0 \le k \le n$,
 the distribution (respectively the generating function) of $Z_n$ given that $Z_k=1$. Therefore, for $k<n$
 \[ f_{k,n} := f_{k+1} \circ \cdots \circ f_n \ .  \]
Our proofs rest on techniques from the paper \cite{Ke}. In particular we use its formula (16), reading
\begin{align} \mathbf P(Z_n>0) =1-f_{0,n}(0)= \Big(\frac 1{ \mu_n} + \sum_{k=1}^n \frac{\varphi_k(f_{k,n}(0))}{\mu_{k-1}}\Big)^{-1} \ , 
\label{probab}
\end{align}
where $\varphi_k(s)$, $0\le s < 1$, is given by the equation $(1-f_k(s))^{-1}= (f_k'(1)(1-s))^{-1} + \varphi_k(s)$. For a detailed discussion  of this {\em shape function} of $f_k$ we refer to \cite{Ke}. We preface the proofs of the theorems with several lemmas.

\paragraph{Lemma 1.} {\em Assume} ($\ast$). {\em Then  condition \eqref{critical} implies
\[ \max_{1\le i\le n}\frac{\nu_i}{\mu_{i-1}} =o(\rho_n) \]
as $n \to \infty$.}

\begin{proof}
From ($\ast$) we have $ \sum_{k > c_{1/2}(1+f_i'(1))} k^2 f_i[k]\le  \sum_{2\le k \le c_{1/2}(1+f_i'(1))} k^2 f_i[k]$, thus
\begin{align}
f_i''(1)\le \sum_{k \ge 2} k^2 f_i[k] \le 2  \sum_{k \le c_{1/2}(1+f_i'(1))} k^2 f_i[k]\le 2c_{1/2} (1+f_i'(1))f_i'(1),
\label{secondfirst}
\end{align}
and therefore for all $i\le n$ and any natural number $k_0$
\begin{align*} \frac{\nu_i}{\mu_{i-1}\rho_n}  &\le 2c_{1/2}  \frac{1+f_i'(1)}{f_i'(1)\mu_{i-1}\rho_n} = 2c_{1/2} \Big(\frac 1{\mu_i\rho_n} + \frac 1{\mu_{i-1} \rho_{n}}\Big) \\
& \le 2c_{1/2}\sup_{k>k_0} \Big(\frac 1{\mu_k\rho_k} + \frac 1{\mu_{k-1} \rho_{k-1}}\Big)+2c_{1/2}\frac 1{\rho_n} \max_{k\le k_0}  \Big(\frac 1{\mu_k} + \frac 1{\mu_{k-1} }\Big)   \ .
\end{align*}
In view of \eqref{critical},  by adapting $k_0$ the first term on the right-hand side can be made arbitrarily small, and similarly the second term by afterwards increasing $n$. This proves the lemma.
\end{proof}

\paragraph{Lemma 2.} {\em Assume} ($\ast$) {\em and let $0< \eta <1$. Then, in case of a critical BPVE  we have for all $0\le i \le n$ with $\rho_i \le  \eta  \rho_n$
\[ \mathbf P(Z_n>0\mid Z_i=1)= \frac {2+o(1)}{\mu_i(\rho_n-\rho_i)}  \]
as $n \to \infty$, where the $o(1)$-term goes to zero uniformly in all $i$ under consideration.}

\begin{proof} Conditioning on the event  $\{Z_i=1\}$ means that we have to replace the terms $\mu_{j}$ for $j \ge i$ in \eqref{probab} by the product $f_{i+1}'(1)\cdots f_{j}'(1)=\mu_j/\mu_i$. Thus, for any $0\le i\le n$
\begin{align}
\mathbf P(Z_n>0\mid Z_i=1) = \Big(\frac {\mu_i}{ \mu_n} + \mu_i\sum_{k=i+1}^n \frac{\varphi_k(f_{k,n}(0))}{\mu_{k-1}}\Big)^{-1} \ .
\label{survivalpr1}
\end{align}

Next we recall from Lemma 8 and formula (7) of \cite{Ke} the approximation
\[ \sum_{k=1}^n \frac{\varphi_k(f_{k,n}(0))}{\mu_{k-1}} = \frac  {\rho_n}2(1+ o(1)) \]
as $n \to \infty$, valid under the assumptions ($\ast$) and \eqref{critical}. Therefore,  because of \eqref{critical}
\begin{align*} \Big| \sum_{k=i+1}^n \frac{\varphi_k(f_{k,n}(0))}{\mu_{k-1}} - \frac{\rho_n-\rho_i}2\Big| &\le \Big| \sum_{k=1}^n \frac{\varphi_k(f_{k,n}(0))}{\mu_{k-1}} - \frac{\rho_n}2\Big|+\Big| \sum_{k=1}^i \frac{\varphi_k(f_{k,n}(0))}{\mu_{k-1}} - \frac{\rho_i}2\Big|\\
&=o(\rho_n) + o(\rho_n) = o(\rho_n)
\end{align*}
uniformly for all $0\le i \le n$. Also $\rho_n =O(\rho_n-\rho_i)$ uniformly for all $i$ satisfying $\rho_i \le  \eta \rho_n$, hence we arrive at
\begin{align}\sum_{k=i+1}^n \frac{\varphi_k(f_{k,n}(0))}{\mu_{k-1}} = \frac{\rho_n-\rho_i}2 (1+o(1)) 
\label{survivalpr2}
\end{align}
uniformly for all $i$ satisfying $\rho_i \le  \eta  \rho_n$.

Also, because of \eqref{critical} we have
\[ \frac 1{\mu_n} = o(\rho_n) = o(\rho_n-\rho_i) \]
uniformly for all $i$ under consideration. Using this formula together with \eqref{survivalpr2} in equation \eqref{survivalpr1} we obtain our claim.
\end{proof}

\paragraph{Lemma 3.} {\em Under the assumptions of Lemma 2 we have
\[ \frac{f_i'(f_{i,n}(0))}{ f_i'(1)} = \exp\Big( - \frac{ \nu_i(2+o(1) )}{\mu_{i-1}(\rho_n-\rho_i)}\Big)  \ ,\]
with the $o(1)$-term applying uniformly to all $0\le i \le n$ fulfilling $\rho_i \le  \eta  \rho_n$.}

\begin{proof}
First let us show that for any numbers $s_{i,n}\in [f_{i,n}(0),1]$ we have
\begin{align}
f_{i,n}''(s_{i,n})=f_{i,n}''(1)(1+o(1))
\label{lemma2a}
\end{align}
as $n\to \infty$, uniformly for all $i$ under consideration. Fix $\varepsilon>0$.
By means of assumption ($\ast$) we obtain with  $s:=s_{i,n}$
\begin{align}
0 &\le f_i''(1)-f_{i,n}''(s) = \sum_{k=2}^\infty k(k-1) f_i[k] (1- s^{k-2})\notag \\
&\le \sum_{k > c_\varepsilon (1+f_i'(1))} k^2 f_i[k] +  \sum_{k \le c_\varepsilon (1+f_i'(1))} k(k-1)f_i[k](k-2)(1-s)\notag\\
&\le  \varepsilon \sum_{k \ge 2} k^2 f_i[k] + c_\varepsilon(1+f_i'(1)) (1-s)\sum_{k \ge 2} k(k-1) f_i[k] \notag\\
&\le \big(2\varepsilon + c_\varepsilon(1+f_i'(1))(1-f_{i,n}(0))\big) f_i''(1) \ .
\label{lemma2b}
\end{align}
Note  that $1-f_{i,n}(0)= \mathbf P(Z_n>0\mid Z_i=1)$. Therefore, from Lemma 2, uniformly for all $i$ under consideration,
\begin{align*}
(1+f_i'(1))(1-f_{i,n}(0)) &= (1+f_i'(1))\frac {2+o(1)}{\mu_i(\rho_n-\rho_i)}\notag \\
&\le \frac{2+o(1)}{1- \eta}\Big( \frac 1{\mu_i\rho_n}+ \frac 1{\mu_{i-1}\rho_n}  \Big) \ .
\end{align*}
As in the proof of Lemma 1, the right-hand side can be made arbitrarily small, uniformly in $i$. This estimate together with \eqref{lemma2b} implies \eqref{lemma2a}.

Now, from \eqref{lemma2a} we get for suitable $s_{i,n}\in [f_{i,n}(0),1]$
\[ f_i'(1)-f_i'(f_{i,n}(0)) = f''(s_{i,n})(1-f_{i,n}(0))= f_i''(1)(1-f_{i,n}(0))(1+o(1)) \]
or
\[ \frac{f_i'(f_{i,n}(0))}{f_i'(1)} = 1 - \frac{f_i''(1)}{f_i'(1)}(1-f_{i,n}(0))(1+o(1)) \ . \]
Noting $f_i''(1)/f_i'(1)= \nu_i\mu_i/\mu_{i-1}$ and  applying Lemma 2  to    $1-f_{i,n}(0)=\mathbf P(Z_n>0\mid Z_i=1)$ we rewrite this formula as
\[ \frac{f_i'(f_{i,n}(0))}{f_i'(1)} = 1 -  \frac{\nu_i(2+o(1)) }{\mu_{i-1}(\rho_n-\rho_i)} \ . \]
Since $\rho_n-\rho_i \ge (1- \eta )\rho_n$, Lemma 1 implies that the right-hand fraction is uniformly of order $o(1)$. Therefore, we may transform the formula into 
\[ \frac{f_i'(f_{i,n}(0))}{f_i'(1)} = \exp\Big( -  \frac{\nu_i(2+o(1)) }{\mu_{i-1}(\rho_n-\rho_i)} \Big)\ , \]
which is our claim.
\end{proof}

\newpage

\begin{proof}[Proof of Theorem 2.]
For the proof of the theorem it is sufficient to show that for any $0< t  < 1$ and for any sequence  of natural numbers $g_n \le n$, $n \ge 1$, fulfilling
\begin{align} \frac{\rho_{g_n}}{\rho_n} \to  t  
\label{theorem2a}
\end{align}
we have
\begin{align} \mathbf P(G_n \le g_n \mid Z_n>0) \to  t   
\label{theorem2b}
\end{align}
as $n \to \infty$.

Let us recall from \cite{Zu} the   formula  for the probability of the event $\{G_n \ge k\}$ with $k<n$, which readily generalizes from GW processes to BPVEs.  (We leave $G_n$ undefined on the event $Z_n=0$, hence $\{G_n \ge k\}\subset \{Z_n >0\}$.)
Note that $Z_k$ takes the value $z>0$ with probability $f_{0,k}[z]$. Then, in order to realize the event $\{G_n \ge k\}$, precisely one of the $z$ individuals in generation $k$ will have offspring in generation $n$.
Therefore
\[ \mathbf P(G_n\ge k\mid Z_k=z)= z \mathbf P(Z_n>0\mid Z_k=1) \mathbf P(Z_n=0 \mid Z_k=1)^{z-1}  \]
and
\begin{align}
\mathbf P(G_n \ge k) &= \sum_{z=0}^\infty z \mathbf P(Z_n>0\mid Z_k=1) \mathbf P(Z_n=0 \mid Z_k=1)^{z-1} f_{0,k}[z] \notag \\
&= \mathbf P(Z_n>0\mid Z_k=1)  f_{0,k}'(\mathbf P(Z_n=0 \mid Z_k=1)) 
\label{proofb}
\end{align}
Moreover, $f_{0,k}'(s)= f_1'(f_{1,k}(s))f_2'(f_{2,k}(s))\cdots f_k'(s)$ and $\mathbf P(Z_n=0 \mid Z_k=1)= f_{k,n}(0)$ for $k \le n$. Therefore, we end up with
\[\mathbf P(G_n \ge k)= \mathbf P(Z_n>0\mid Z_k=1)  \prod_{i=1}^k f_i'(f_{i,n}(0)) \ . \]

We use this formula with $k$ replaced by $g_n$. Since these numbers satisfy \eqref{theorem2a}, we have $\rho_{g_n}\le \eta \rho_g$ with some $\eta \in (t,1)$ and $n$ sufficiently large. Thus we may use the approximations given by Lemmas 2 and 3 uniformly for all $i \le g_n$ yielding
\[ \mathbf P(G_n \ge g_n)= \frac{2+o(1)}{\mu_{g_n}(\rho_n-\rho_{g_n})} \exp\Big( - \sum_{i=1}^{g_n} \frac{ \nu_i(2+o(1))}{\mu_{i-1}(\rho_n-\rho_i)} \Big)\prod_{i=1}^{g_n} f_i'(1)\ . \]
Note that the term $\mu_{g_n}$ and the product cancel out. Due to the uniform convergence of the $o(1)$-terms inside of the sum, we may replace the factors $2+o(1)$ (which depend on $i$) by one factor $2+o(1)$ outside of the sum. Note also from \eqref{murho} that  $\nu_i/\mu_{i-1}= \rho_i-\rho_{i-1}$ implying
\[ \sum_{i=1}^{g_n} \frac{ \nu_i(2+o(1))}{\mu_{i-1}(\rho_n-\rho_i)} = (2+o(1)) \sum_{i=1}^{g_n} \frac 1{1- \rho_i/\rho_n} \Big(\frac{\rho_i}{\rho_n} -\frac{\rho_{i-1}}{\rho_n}\Big)\ . \]
The   right-hand sum may be viewed as a Riemann approximation to  $\int_0^{\rho_{g_n}/\rho_n}(1-x)^{-1}dx$, using the partition $\rho_i/\rho_n$, $i=0,\ldots,g_n$. Lemma 1 shows, that  its mesh converges to zero. Also taking  \eqref{theorem2a} into account we obtain
\[\sum_{i=1}^{g_n} \frac{ \nu_i(2+o(1))}{\mu_{i-1}(\rho_n-\rho_i)} = (2+ o(1))\int_0^ t  \frac {dx}{1-x}=-(2+o(1))\log(1- t ) \ , \]
consequently
\[ \mathbf P(G_n \ge g_n) =  \frac{2+o(1)}{\rho_n-\rho_{g_n}} \exp \big((2+o(1))\log(1- t )\big)= \frac{2+o(1)}{\rho_n-\rho_{g_n}}(1- t )^2 \ .\]
Again using Lemma 2 with $i=0$ and \eqref{theorem2a}  we arrive at
\[\mathbf P(G_n \ge g_n\mid Z_n>0) =  \frac {\rho_n(1+o(1))}{\rho_n-\rho_{g_n}} (1- t )^2 = 1- t  + o(1) \ . \]
Finally, because of Lemma 1 and \eqref{theorem2a} we have $\rho_{g_n+1}/\rho_n= t  +o(1)$, too, therefore we may replace $g_n$ by $g_n+1$ in our considerations. Hence we obtain as well
\[\mathbf P(G_n > g_n\mid Z_n>0)= 1- t  + o(1),\]
which implies the desired assertion \eqref{theorem2b}.
\end{proof}

\begin{proof}[Proof of Theorem 3] Fix some $\eta \in (0,1)$.
Let $0=T_{0,n}<T_{1,n} < \cdots $ denote the jump times of the process $X_n(t):= Z_{k_n(t),n}$, and $0=T_0<T_1 < \cdots$ the jump times of the limiting process $X(t):=Y(\log \tfrac 1{1-t})$, $0\le t < 1$.  Since otherwise the processes have constant paths, and since $X(t)$ possesses only finitely many jumps during the restricted time interval $[0,\eta]$, it is sufficient for the proof to show that  for all $j \ge 1$ the random variables $X_n(T_{j,n})\cdot 1_{T_{j,n} \le \eta}$ converge to $X(T_j)=j+1$ in probability and that the random vectors $(T_{1,n}\cdot 1_{T_{0,n} \le \eta},  \ldots , T_{j,n}\cdot 1_{T_{j-1,n} \le \eta})$ converge in  distribution to the corresponding limiting vector as $n \to \infty$. To see this, one has to show convergence of the finite-dimensional distributions and  relative compactness of the sequence $(X_n)_{n \ge 1}$. The former requirement follows from 
\[ \{ T_{j,n} >t, X_n(T_{j,n})=j+1\} \subset \{ X_n(t)\le j+1\} \subset  \{ T_{j,n} >t\} \ , \]
and the latter one e.g. by means of \cite[Theorem 7.2]{EtKu}. Namely, the compact containment condition arises from the monotonicity of the sample paths of $X_n$, and for the modulus of continuity $w'$ we have $w'(X_n, \delta, \eta)=0$ for $0<\delta < \min\{ T_{j+1,n}\wedge \eta-T_{j,n} :  j \ge 0 , T_{j,n} < \eta\}$.

The standard Yule process is a pure birth process characterized by the property that the waiting time at any state $j \in \mathbb N$, given the previous waiting times,  is equal in distribution to the minimum of $j$ independent, standard exponential random variables. Hence, this property holds for the terms \[\log (1-T_j)^{-1}-\log (1-T_{j-1})^{-1}= \log \Big(1- \frac{T_j-T_{j-1}}{1-T_{j-1}}\Big)^{-1} . \]
This readily translates to the characterization of the process $X(t)$, $0\le t < 1$, as a (non-homogeneous) pure birth process with the property that for any $j \ge 1$, given the values of $T_1, \ldots , T_{j-1}$, the random variable   $(T_j-T_{j-1})/(1-T_{j-1})$ has the distribution of the minimum of $j$ independent random variables, each   uniformly distributed on the interval $(0,1)$. We are going to check this property for the limiting process. 

To ease  notation we set $\mathbb P_n(\cdot)= \mathbf P(\cdot \mid Z_n>0)$. We prepare the proof by two approximations. Let $0\le t < t'\le \eta$. Then,
\begin{align*}
\mathbb P_n( T_{1,n}>t' \mid X_n(t)=1) = \frac {\mathbb P_n(X_n(t')=1)}{\mathbb P_n(X_n(t)=1)} = \frac{\mathbf P(G_n \ge k_n(t'))}{\mathbf P(G_n \ge k_n(t))} \ .
\end{align*}
Using Theorem 2 it follows that in the limit $n \to \infty$
\[\mathbb P_n( T_{1,n}>t' \mid X_n(t)=1) \sim \frac{1- \rho_{k_n(t')}/\rho_n}{1- \rho_{k_n(t)}/\rho_n} \ , \]
and because of  $\rho_{k_n(u)} \le u \rho_n \le \rho_{k_n(u)+1}$ with $u=t,t'$  we obtain by means of  Lemma 1 
\begin{align} 
\mathbb P_n( T_{1,n}>t' \mid X_n(t)=1) \to \frac{1-t'}{1-t} 
\label{proofa}
\end{align}
uniformly in all $t \le \eta$. 

Similarly, 
\begin{align*}
\mathbb P_n(X_n( T_{1,n})=2 \mid X_n(t)=1,  T_{1,n}=t')=  \frac{\mathbf P(G_n = k_n(t')-1,Z_{k_n(t'),n}=2)}{\mathbf P(G_n = k_n(t')-1)} \ .
\end{align*}
Let us write $k'=k_n(t')$ for short. In order to realize the event $\{G_n =  k' -1,Z_{ k' ,n}=2\}$ it is required that from  $z$ individuals in generation $ k' -1$  just one  has descendants in generation $n$, and from its $z'$-many offspring  precisely two succeed in that respect. Given $z$ and $z'$ this event has probability
\begin{align*} 
z f_{0, k' -1}[z] \mathbf P(Z_n=0&\mid Z_{ k' -1}=1)^{z-1} \\&\cdot \frac{z'(z'-1)}2 f_{ k' }[z'] \mathbf P(Z_n>0\mid Z_{ k' }=1)^2\mathbf P(Z_n=0\mid Z_{ k' }=1)^{z'-2} \ .
\end{align*}
Summing over $z$ and $z'$ we get
\begin{align*}\mathbf P(G_n =  k' -1,Z_{ k' }=2)= \frac 12 \mathbf P(Z_n>0\mid Z_{ k' }=1)^2 f_{0, k'-1 }'\big(f_{k'-1,n}(0)\big)f_{ k' }''\big(f_{k',n}(0)\big)\ .
\end{align*}
Moreover, from \eqref{proofb}
\begin{align*} 
\mathbf P(G_n=k'-1)&= \mathbf P(G_n \ge k'-1)- \mathbf P(G_n \ge k') \\
&=(1-f_{k'-1,n}(0))f_{0,k'-1}'(f_{k'-1,n}(0)) - (1-f_{k',n}(0))f_{0,k'}'(f_{k',n}(0)) \ ,
\end{align*}
and because of $f_{0,k'}'(f_{k',n}(0))= f_{0,k'-1}'(f_{k'-1,n}(0))f_{k'}'(f_{k',n}(0))$ and $f_{k'-1,n}(0)=f_{k'}(f_{k',n}(0))$
\begin{align*}
\mathbf P(G_n=k'-1)&=f_{0,k'-1}'(f_{k'-1,n}(0)) \cdot\big(1- f_{k'}(f_{k',n}(0))- (1-f_{k',n}(0))f_{k'}'(f_{k',n}(0))\big)\\
&= f_{0,k'-1}'(f_{k'-1,n}(0))\cdot \frac 12 f_{k'}''(s_{t',n}) (1-f_{k',n}(0))^2 \ ,
\end{align*}
with suitable $ s_{t',n} \in [f_{k',n}(0), 1]$. Consequently,
\[ \frac{\mathbf P(G_n = k_n(t')-1,Z_{k_n(t'),n}=2)}{\mathbf P(G_n = k_n(t')-1)}= \frac{f_{ k_n(t') }''\big(f_{k_n(t'),n}(0))}{f_{k_n(t')}''(s_{t',n})} \ . \]
Since $\rho_{k_n(t')} \le t'\rho_n\le \eta \rho_n$, we may use \eqref{lemma2a} implying that
\begin{align} \mathbb P_n(X_n( T_{1,n})=2 \mid X_n(t)=1,  T_{1,n}=t') =1 + o(1)  
\label{proofc}
\end{align}
uniformly in all $0 \le t < t' \le \eta$.

We are now ready to prove that for any $j\ge 1$ the random vector
\begin{align}( T_{1,n}\cdot 1_{T_{0,n}\le\eta}, X_n(T_{1,n})\cdot 1_{T_{1,n}\le\eta}, \ldots, T_{j,n}\cdot 1_{T_{j-1,n}\le \eta}, X_n(T_{j,n})\cdot 1_{T_{j,n}\le\eta}) 
\label{proofd}
\end{align}
has the announced limiting behaviour. We proceed by induction on $j$. For $j=1$ the claim follows immediately from \eqref{proofa} and \eqref{proofc} by setting $t=0$ (taking into account the uniform convergence in \eqref{proofc} concerning $t'$). For the induction step suppose that our claim is valid for $j-1$. Then the induction hypothesis can be applied to the left part in \eqref{proofd}. In particular,  with increasing $n$ the probability approaches 1 that the term $X_n(T_{j-1,n})\cdot 1_{T_{j-1,n} \le \eta}$ takes the value $j$. Then precisely $j$ individuals in generation $k_n(T_{j-1,n})$ will have descendants until generation $n$. Due to the  properties of the underlying branching process the induced $j$ family trees are independent and identically distributed. Therefore, the waiting times $W_{1,j,n}$, \ldots, $W_{j,j,n}$ within these subtrees up to their first branching events are i.i.d. random variables. Due to the Markov property of the reduced process 
their distributions, given $T_{j-1,n}=t$, depend only on the distributions $f_k$ with $k_n(t)\le k \le n$, and not on $j$. In particular, for $1 \le i \le j$ and $0\le t < t' \le 1$
\[ \mathbb P_n( W_{i,j,n}>t' -t\mid T_{j-1,n}=t) = \mathbb P_n( T_{1,n}>t'\mid X_n(t)=1 ) \ ,\]
so that because of \eqref{proofa} the conditional distribution of $W_{i,j,n}$ is asymptotically uniform on $[0,1-t]$, and the asymptotic distribution of $T_{j,n}-T_{j-1,n}= \min\{W_{1,j,n}$, \ldots, $W_{j,j,n}\}$ given $T_{j-1,n}$ coincides in distribution with the minimum of $j$ uniform random variables on $[0,1-T_{j-1,n}]$. 
  Also, with probability going to 1 all $W_{1,j,n}$, \ldots, $W_{j,j,n}$  differ from each other. Then, in view of \eqref{proofc} and in case of $T_{j,n}\le \eta$,  the process $X_n$ will increase at time $T_{j,n}$  just by 1 and $X_n(T_{j,n})$ will take the value $j+1$ in the limit $n\to \infty$. This finishes the induction and concludes the theorem's proof.
\end{proof}

\paragraph{Acknowledgement.} It is my pleasure to dedicate this work to Vladimir Vatutin and Andrey Zubkov.

\end{document}